\documentclass[a4paper]{article}
\usepackage[left=3.0cm,right=2.0cm,top=2.5cm,bottom=2.5cm]{geometry}
\usepackage{xcolor}
\usepackage{graphicx}
\usepackage{cite}
\usepackage{amsmath,amssymb,amscd,amsxtra,amsfonts}

\newtheorem{example}{Example}
\providecommand{\Div}{\operatorname{div}}
\providecommand{\curl}{\operatorname{{\bf curl}}}
\newcommand{\Bx}{{\boldsymbol{x}}}
\newcommand{\Bv}{{\boldsymbol{v}}}
\newcommand{\Ba}{{\boldsymbol{a}}}
\newcommand{\D}{\mathrm{d}}
\newcommand{\varphibf}{\boldsymbol{\varphi}}

\newcommand{\Vb}{{\mathbf{b}}}

\newcommand{\Ve}{{\mathbf{e}}}

\newcommand{\Vg}{{\mathbf{g}}}

\newcommand{\Vn}{{\mathbf{n}}}

\newcommand{\Vr}{{\mathbf{r}}}

\newcommand{\Vu}{{\mathbf{u}}}

\newcommand{\Vw}{{\mathbf{w}}}
\newcommand{\Vx}{{\mathbf{x}}}
\newcommand{\Vy}{{\mathbf{y}}}

\newcommand{\VA}{{\mathbf{A}}}
\newcommand{\VB}{{\mathbf{B}}}

\newcommand{\VE}{{\mathbf{E}}}

\newcommand{\VH}{{\mathbf{H}}}

\newcommand{\VJ}{{\mathbf{J}}}

\newcommand{\BB}{{\boldsymbol{B}}}

\newcommand{\BF}{{\boldsymbol{F}}}

\newcommand{\BH}{{\boldsymbol{H}}}
\newcommand{\BI}{{\boldsymbol{I}}}

\newcommand{\BK}{{\boldsymbol{K}}}
\newcommand{\BL}{{\boldsymbol{L}}}
\newcommand{\BM}{{\boldsymbol{M}}}

\newcommand{\BP}{{\boldsymbol{P}}}

\newcommand{\BX}{{\boldsymbol{X}}}

\newcommand{\Ct}{\mathcal{T}}
\newcommand{\bbN}{\mathbb{N}}
\newcommand{\bbZ}{\mathbb{Z}}
\newcommand{\bbT}{\mathbb{T}}
\newcommand{\bbV}{\mathbb{V}}
\newcommand{\bbS}{\mathbb{S}}
\newcommand{\bbR}{\mathbb{R}}

\newcommand*{\N}[1]{\left\|{#1}\right\|}

\newcommand*{\Lp}[2][\defaultdomain]{L^{#2}({#1})}
\newcommand*{\Lpv}[2][\defaultdomain]{\BL^{#2}({#1})}
\newcommand*{\NLp}[3][\defaultdomain]{\N{#2}_{\Lp[#1]{#3}}}

\newcommand*{\Ltwov}[1][\defaultdomain]{\Lpv[#1]{2}}
\newcommand*{\NLtwo}[2][\defaultdomain]{\NLp[#1]{#2}{2}}

\newcommand*{\Hcurl}[1][\defaultdomain]{\boldsymbol{H}(\curl,{#1})}
\newcommand*{\bHcurl}[2][\defaultdomain]{\boldsymbol{H}_{#2}(\curl,{#1})}
\newcommand*{\zbHcurl}[1][\defaultdomain]{\bHcurl[#1]{0}}
\newcommand*{\NHcurl}[2][\defaultdomain]{\N{#2}_{\Hcurl[#1]}}

\title{A Conservative Finite Element Solver for MHD Kinematics equations: Vector Potential method and Constraint Preconditioning}
\author{Xujing Li
\thanks{School of Mathematics, Hangzhou Normal University, Hangzhou 311121, China. (lixujing@hznu.edu.cn)}
\and Lingxiao Li
\thanks{\textbf{Corresponding author}: Institute of Applied Physics and Computational Mathematics (IAPCM), Fenghaodong Road, Haidian District,
Beijing 100094, China. (lilingxiao@lsec.cc.ac.cn)}}
\date{2021}
\begin{document}
\maketitle
\begin{abstract}
  A new conservative finite element solver for the three-dimensional steady magnetohydrodynamic (MHD)
  kinematics equations is presented.
  The solver utilizes magnetic vector potential and current density as solution variables,
  which are discretized by H(curl)-conforming edge-element and H(div)-conforming face element respectively.
  As a result, the divergence-free constraints of discrete current density and magnetic induction
  are both satisfied. Moreover the solutions also preserve the total magnetic helicity.
  The generated linear algebraic equation is a typical dual saddle-point problem that is ill-conditioned and indefinite.
  To efficiently solve it, we develop a block preconditioner based on constraint preconditioning framework
  and devise a preconditioned FGMRES solver. Numerical experiments verify the conservative properties,
  the convergence rate of the discrete solutions and the robustness of the preconditioner.
\end{abstract}
\textbf{Key Workds:}
MHD; Divergence-free conditions; Block preconditioner; Constraint preconditioning; Magnetic helicity conservation.
\section{Introduction}
\label{sec:introduction}
Magnetohydrodynamics (MHD) has broad applications in our real world.
It describes the interaction between electrically conducting fluids
and magnetic fields, which plays an important role in magnetic confined fusion \cite{jardin2010},
Z-pinch \cite{pereira2000},
astrophysics and liquid metals \cite{dav01}.
In this paper, we are studying the conservative finite element method and efficient iterative solver
for the following steady MHD kinematics equations
\begin{subequations}\label{eq:mhd}
\begin{align}
\curl\VE                  &=   \mathbf{0}\qquad\hbox{in}\;\;\Omega,   \label{eq:Faraday}\\
\curl\VH - \VJ - \VJ_s    &=\mathbf{0}   \qquad\hbox{in}\;\;\Omega,   \label{eq:ns}\\
\sigma(\VE+\Vw\times\VB)  &= \VJ         \qquad\hbox{in}\;\;\Omega,   \label{eq:OhmLaw}\\
\Div\VJ =0,   \quad\Div\VB&=0            \qquad \hbox{in}\;\;\Omega.  \label{eq:div}
\end{align}
\end{subequations}
where $\VE$ and $\VH$ are the electric field and the magnetic field respectively. $\VB$ is the magnetic flux density.
$\Vw$ is the prescribed velocity field, $\VJ, \VJ_s$ the induced current and source current.
We assume that $\Omega$ is a bounded, simply-connected, and Lipschitz polyhedral domain with boundary $\Gamma=\partial\Omega$.
When $\Vw \equiv 0$ in the domain, this model reduces to the classical eddy currents model \cite{buff2000, hip02}.
The equations in \eqref{eq:mhd} are complemented with the following constitutive equation
\begin{equation}
\VB=\mu\VH
\end{equation}
The MHD kinematics equations have interest applications in the field of dynamo theory \cite{liu2007, mini2005, phil2014}.
Such applications contain MHD generators, dynamo of the sun, brine and the geodynamo.
Combined with the momentum equations MHD kinematics equation becomes the full MHD equation,
so its efficient solver constitutes a core part of the MHD solver.
In the present work we propose a new finite element method which preserves the
divergence-free conditions for both magnetic induction and current density at the same time. Moreover,
a robust block preconditioner
is developed based on constraint preconditioning framework \cite{KGW2000, C2002}.

There already exists extensive papers in the literature to study numerical methods for MHD.
Now we give a short review but not complete reference list for relevant topics.
In \cite{gun91}, Gunzburger et al studied well-posedness and the finite element
method for the stationary incompressible MHD
equations. The magnetic field is discretized by the $H^1(\Omega)$-conforming
finite element method. In \cite{G2000NM}, Gerbeau introduced a stabilized finite element method for the incompressible MHD.
We also refer to \cite{ger06} for a systematic
analysis on finite element methods for incompressible MHD
equations. In 2004, Sch\"{o}tzau \cite{sch04}
proposed a novel mixed finite element method to solve the stationary
incompressible MHD equations where edge elements are used to solve the
magnetic field.
In 2010, Greif et al. extended the work in \cite{sch04} by H(div)-conforming face elements for velocity such
that $\Div\Vu_h = 0$ holds exactly \cite{gre10}.
Here we represent the magnetic induction by vector potential variable such that $\VB=\curl\VA$ and use edge element
to discretize $\VA$. As a result in the discrete level $\Div\VB_h = 0$ is naturally satisfied.
The theoretical foundation can be found in previous publication such as \cite{ABDG1998}. For error analysis of finite element method,
we refer to \cite{H2015IMA} for Euler semi-implicit scheme and \cite{SMF2019} for penalty-based finite element methods.
In \cite{ABGM2021}, Alvarez, Bokil, Gyrya and Manzini devised a novel virtual element method
for time-dependent MHD kinematics equations which is similar to
the physical model considered in our present work.
Moreover, Stasyzyn and Elstner in \cite{SEJCP2015} introduced a smoothed
particle magnetohydrodynamics algorithms, where
magnetic vector potential $\VA$ with Coulomb gauge is implemented.

In recent years, exactly divergence-free approximations for $\VJ$ and $\VB$ have attracted
more and more interest in numerical simulation. For the current density $\VJ$ we would like to mention the
current density-conservative finite volume methods of Ni et al. for the inductionless MHD model on both structured and unstructured
grids \cite{ni12, ni07-1, ni07-2}. In these work, the authors showed that when the applied magnetic field is constant,
the discrete Lorentz force in the momentum equation can precisely conserve the total momentum when the
current density is divergence-free. And they suggested that only the divergence-free schemes which conserve the total momentum in the discrete
level can obtain accurate result for MHD flow at large Hartmann numbers.
In \cite{L2019SIAM}, Li et al developed an charge-conservative finite element method for inductionless MHD equations.
In fact charge-conservative property is an important constraint in plasma physics,
thus for accurate numerical simulation, the discrete methods should preserve this feature.
In the present work, H(div)-conforming element is used for discrete current density $\VJ_h$
to reach the goal $\Div\VJ_h=0$, which is the same as \cite{L2019SIAM}.

The importance of divergence-free condition for $\VB_h$ has been discussed for a long period. From Ramshaw \cite{R1983JCP}, Evans \cite{AH1988}
to T\'{o}th \cite{T2000JCP} one can see thorough arguments for this property.
For this point, we would like to mention the pioneering work in
\cite{RRW2005,RWWS2007}. In \cite{RRW2005}, using edge element for $\VE$ and H(div)-conforming element for $\VB$,
Rieben et al. developed a high order finite element solver for time-dependent Maxwell equations, where the discrete magnetic induction
is exactly divergence-free. Motivated by the concept of
differential form \cite{hip02}, then in \cite{RWWS2007}, Rieben,
White, Wallin and Solberg of LLNL successfully
extended the ideas of \cite{RRW2005} to 3D compressible MHD equations in the ALE framework. Again they achieved the precise divergence-free
conditions for $\VB_h$.
For incompressible MHD equations, Hu et al. in \cite{hu17} discretize the electric field $\VE$ by edge elements
and the magnetic induction $\VB$ by H(div)-conforming elements such that $\Div\VB_h = 0$ is achieved.
In \cite{hi2018}, for time-dependent MHD equations, Hiptmair et al. use temporal gauge to represent the electric field by $\VE=-\partial_t\VA$ and
magnetic induction by $\VB = \curl\VA$. With edge element for $\VA_h$,  the Gauss's law for $\VB_h$ is satisfied and they also proved
the convergence of the finite element solutions.
In \cite{ABGM2021}, Alvarez et al. presented a novel virtual element method for resistive
MHD where the divergence of $\VB_h$ is automatically zero.
Very recently, Li et al. \cite{lxliPHD,LZZ2021} proposed a constrained transport divergence-free
finite element method for incompressible MHD equations, where the authors achieve the conditions $\Div\VJ_h = \Div\VB_h = 0$ at the same time using magnetic field $\VH$ and vector potential
$\VA$ as variables. Different from the work in \cite{LZZ2021}, in this paper,
we use current density $\VJ$ and $\VA$ as main variables
and the methods in \cite{L2019SIAM} is incorporated
to realize the conditions $\Div\VJ_h = 0$.

Another objective of this paper is to propose a preconditioned iterative
method to solve the algebraic systems associated with the proposed divergence-free finite element solver.
In this procedure the key ingredient is efficient preconditioning \cite{ESW2014}. For MHD equations, large number of
studies exist in the literature, such as
\cite{phil2014,PSCEP2016,C2008POP} and references therein, on block preconditioners using approximate Schur complements techniques.
We also refer to the work in \cite{shadid2010,shadid2016} for algebraic multigrid methods and in \cite{adler2016} for geometric multigrid method.
In particular we point out that in \cite{phil2014},
Phillips and Elman constructed an efficient block preconditioner for steady MHD kinematics equations
with $\VB$ and an extra multiplier $r$ as variables which is a sub-block of the model in \cite{sch04}.
In the present paper, we will derive a
block preconditioner based on the constraint preconditioning framework \cite{KGW2000, C2002},
which is different from the techniques mentioned above.

The paper is organized as follows:
In section 2, we introduce the dimensionless model using magnetic vector potential $\VA$
and electrical potential $\phi$.
In section 3, we introduce a variational formulation for the MHD kinematics equations and show that the discrete
formulation can preserve the divergence-free properties for $\VB_h$ and $\VJ_h$ precisely.
Besides the magnetic helicity is also preserved.
In section 4, from the matrix level, we introduce the constraint preconditioning framework from \cite{KGW2000, C2002}
for our dual-saddle problem and give some eigenvalues discussions.
A block preconditioner is developed in this section.
In section 5, numerical experiments are conducted to verify the conservation of the discrete solutions,
the convergence rate of the finite element solver,
and to demonstrate the optimality and the robustness of the iterative solver.
In section 6, some conclusions and further investigations are pointed out.

Throughout the paper we denote vector-valued quantities by boldface notation, such as $\Ltwov:=(L^2(\Omega))^3$. In the following, we assume the
physical parameters $\sigma,\mu$ are constants despite that the solver can be adapted to variable coefficient case.

\section{A dimensionless vector potential formulation}
In this section we will derive the vector potential formulation for our finite element iterative solver.
First note that $\curl\VE = \mathbf{0}$ in \eqref{eq:Faraday} so we have
\[\VE = - \nabla\phi\]
where $\phi$ is generally called electric potential.
Due to the theory in \cite{ABDG1998, LZZ2021} we can represent the magnetic induction by vector potential $\VA$
such that
\[\VB = \curl\VA,\quad \Div\VA = 0\]
where the second divergence constraint for $\VA$ is called Coulomb's gauge condition.
In short we have
\begin{equation}\label{eq:trans}
\VE = - \nabla\phi,\quad \VB = \curl\VA,\quad \Div\VA = 0
\end{equation}

Using the transformation \eqref{eq:trans} and the generalized Ohm's law $\VJ = \sigma(\VE+\Vw\times\VB)$,
one will obtain the following new formulation
\begin{subequations}\label{eq:new-mhd}
\begin{align}
\sigma^{-1}\VJ + \nabla\phi - \Vw\times\curl\VA  = \mathbf{0},&\quad \mathrm{in}~~\Omega
\label{current}\\
\Div\VJ = 0, &\quad \mathrm{in}~~\Omega
\label{conserv}\\
-\VJ + \curl\mu^{-1}\curl\VA = \VJ_s,&\quad \mathrm{in}~~\Omega
\label{vectorp}\\
\Div\VA = 0, &\quad \mathrm{in}~~\Omega
\label{aguage}
\end{align}
\end{subequations}
Let $L$, $u_0$, $B_0$ and $\sigma_0$ be the characteristic length, characteristic velocity, characteristic magnetic flux density and reference conductivity respectively and make the following scaling
\begin{equation}\label{mhd-scal}
\Bx\leftarrow \Bx\frac{1}{L},~\VJ\leftarrow\VJ\frac{1}{\sigma_0B_0 u_0},
~\Vg\leftarrow\VJ_s\frac{1}{\sigma_0B_0 u_0},
~\phi\leftarrow \phi\frac{1}{B_0u_0 L},
~\sigma\leftarrow \sigma\frac{1}{\sigma_0}
\end{equation}
Then we can get the desired dimensionless formulation
\begin{subequations}\label{eq:steady-mhd}
\begin{align}
\sigma^{-1}\VJ + \nabla\phi - \Vw\times\curl\VA  = \mathbf{0},&\quad \mathrm{in}~~\Omega\label{eq:steady-current}\\
\Div\VJ = 0,                                                  &\quad \mathrm{in}~~\Omega\\
-\VJ + \textsf{Rm}^{-1}\curl\curl\VA = \Vg,                   &\quad \mathrm{in}~~\Omega\label{eq:steady-A}\\
\Div\VA = 0,                                                  &\quad \mathrm{in}~~\Omega
\end{align}
\end{subequations}
where $\textsf{Rm} = \mu\sigma_0 L u_0$ is the magnetic Reynolds number.
For simplicity system of equation \eqref{eq:steady-mhd} is complemented with the following homogeneous boundary conditions
\begin{equation}\label{eq:bc}
\phi = 0,\quad \VA\times\Vn = \mathbf{0} \qquad\hbox{on}\;\;\Gamma.
\end{equation}
Since $\VE = -\nabla\phi$ and $\VB = \curl\VA$, from \eqref{eq:bc} one can deduce that
(see Theorem 3.17 of \cite{ABDG1998} and Section 3.7 of \cite{monkbook})
\begin{equation}
\VE\times\Vn = 0,\quad \VB\cdot\Vn = 0 \qquad\hbox{on}\;\;\Gamma.
\end{equation}
which is precisely the perfect conducting boundary condition \cite{monkbook}.

In the following, we will focus on the steady system \eqref{eq:steady-mhd} and devise a conservative finite element method.
More importantly, we will construct a block preconditioner from the constraint preconditioning framework.
For convenience electrical resistivity $\eta$ will be used instead of $\sigma^{-1}$ and $\nu_m$ instead of $\textsf{Rm}^{-1}$ in some places.

\section{Variational formulation and mixed finite element method}
First we will introduce the Hilbert spaces and Sobolev norms used in
this paper. Let $L^2(\Omega)$ be the usual Hilbert space of square
integrable functions equipped with the following inner product and
norm:
\begin{eqnarray}
(u,v):=\int_{\Omega}u(\Bx)\,v(\Bx)\D\Bx \quad \hbox{and} \quad
\|u\|_L^2(\Omega):=(u,u)^{1/2}. \nonumber
\end{eqnarray}

Define $H^m(\Omega):=\{v\in L^2(\Omega): D^{\alpha}v\in
L^2(\Omega),|\alpha|\le m\}$ where $\alpha$ represents non-negative triple
index. Let $H^1_0(\Omega)$ be the subspace of $H^1(\Omega)$ whose
functions have zero traces on $\Gamma$. We define the spaces of functions having square integrable curl by
\begin{eqnarray}
\Hcurl&:=&\{\Bv\in\Ltwov\,:\;\curl\Bv\in \Ltwov\}, \nonumber\\
\zbHcurl&:=&\{\Bv\in\Hcurl\,:\;\Vn\times\Bv=0\;\;
\hbox{on}\;\Gamma\}, \nonumber
\end{eqnarray}
which are equipped with the following inner product and norm
\begin{equation*}
(\Bv,\Vw)_{\Hcurl}:=(\Bv,\Vw)+(\curl\Bv,\curl\Vw), \;\;
\NHcurl{\Bv}:=\sqrt{(\Bv,\Bv)_{\Hcurl}}\;.
\end{equation*}
here $\Vn$ denotes the unit outer normal to $\Gamma$. We also use the usual Hilbert space $\BH(\Div,\Omega)$ indicating square integrable divergence.
We will introduce the notations $V^p$ $(p=0,1,2,3)$ for the Hilbert spaces mentioned above as following:
\[V^0=H_0^1(\Omega),~V^1=\zbHcurl,~V^2=\BH(\Div,\Omega),~V^3=L^2(\Omega)\]
where the superscript $p$ of $V^p$ indicates that the corresponding physical field is
\textsf{differential $p$-form} \cite{hip02,arnold2020}.
It is known that the spaces pair $(V^2,V^3)$ and $(V^1,V^0)$ both satisfy the LBB conditions \cite{BrennerFEM,monkbook}.

It is a standard procedure in variational theorem to obtain the continuous mixed variational formulation for the system \eqref{eq:steady-mhd}.
An extra Lagrange multiplier $r \in H_0^1(\Omega)$ for magnetic vector potential
is introduced for the sake of well posedness.
This skill has been used  for $\VB$ variable to obtain a well-posed variational form in \cite{sch04}.
The proposed continuous
weak form reads as:

Find $(\VJ, \phi, \VA, r) \in V^2 \times V^3 \times V^1 \times V^0$ such that the following weak formulation holds
\begin{subequations}\label{weaka}
\begin{align}
\eta(\VJ,\varphibf) - (\phi,\Div\varphibf) - (\Vw\times\curl\VA, \varphibf)   &= 0,\label{weaka:J}\\
-(\Div\VJ, \psi) &= 0, \label{weaka:Phi} \\
-(\VJ, \Ba) + \nu_m(\curl\VA, \curl\Ba) + (\nabla r, \Ba) &= (\Vg, \Ba), \label{weaka:A}\\
(\VA, \nabla s) &= 0. \label{weaka:r}
\end{align}
\end{subequations}
for any $(\varphibf, \psi, \Ba, s) \in V^2 \times V^3 \times V^1 \times V^0$.

Based on \eqref{weaka}, we define the bilinear forms by
\[a_1(\VJ,\varphibf)=\eta\int_{\Omega} \VJ\cdot\varphibf,\qquad a_2(\VA,\Ba) = \nu_m\int_{\Omega}\curl\VA\cdot\curl\Ba\]
\[d_1(\VJ,\phi)=-\int_{\Omega}\phi\Div\VJ,\qquad d_2(\VA,r)=\int_{\Omega} \VA\cdot\nabla r\]
and the trilinear form by
\[c(\Vw;\VA,\varphibf)= - \int_{\Omega} \Vw\times\curl\VA\cdot\varphibf\]
Note if we let $a = \nabla r$ in \eqref{weaka:A}, considering $r\in H_0^1(\Omega)$ and $\Div\VJ = 0$, we will have
\[(\nabla r, \nabla r) = (\Vg, \nabla r) = -(\Div\Vg, r)\]
which means $r=0$ in the domain if the divergence of the source term $\Vg$ vanishes.
Generally in practical application the source current density $\Vg$ is divergence-free,
thus Lagrangian multiplier $r$ does not influence the exact solutions as in \cite{sch04}.

\subsection{Mixed finite element method}
Let $\mathcal{T}_h$ be a shape-regular tetrahedral triangulation of $\Omega$, with $h$ the grid size if the partition is quasi-uniform.
We will use finite element spaces which are all conforming, namely
\begin{equation*}
V_h^2 \subseteq V^2,~~V_h^3 \subseteq V^3,~~V_h^1 \subseteq V^1,~~V_h^0 \subseteq V^0
\end{equation*}
Moreover we need the finite element space pair $(V_h^2,V_h^3)$ and $(V_h^1,V_h^0)$ both
satisfy the inf-sup conditions \cite{BrennerFEM,monkbook}.
For simplicity, in our numerical experiments, for $V_h^2$ we use the $\BH(\Div,\Omega)$-conforming piecewise linear finite element \cite{xin2013}
\[V_h^2 = \{\varphibf_h\in \BH(\Div,\Omega): \varphibf_h|_K\in\BP_1(K), ~ K\in \mathcal{T}_h\}\]
For $V_h^3$ we use the piecewise constants finite element
\[V_h^3 = \{\psi_h\in L^2(\Omega): \psi_h|_K\in P_0(K), ~ K\in \mathcal{T}_h\}\]
The finite element for $\VA$ is the first order N\'{e}d\'{e}lec edge element space \cite{ne86}
\[V_h^1 = \{\Ba_h\in \zbHcurl: \Ba_h|_K\in\BP_1(K),~ K\in \Ct_h\}\]
The Lagrangian finite element space $V_h^0$ for $r_h$ is defined by
\[V_h^0 = \{s_h \in H_0^1(\Omega): s_h|_K \in P_2(K),~K \in \Ct_h\}\]

The discrete mixed finite element scheme to solve the continuous formulation \eqref{weaka}
reads as:

Find $(\VJ_h, \phi_h, \VA_h, r_h) \in V_h^2 \times V_h^3 \times V_h^1 \times V_h^0$ such that the following weak formulation holds
\begin{subequations}\label{full-dis}
\begin{align}
a_1(\VJ_h,\varphibf_h)+ d_1(\varphibf_h, \phi_h) + c(\Vw;\VA_h, \varphibf_h) &= 0,\label{full-dis:J}\\
d_1(\VJ_h, \psi_h) &= 0,\label{full-dis:Phi}\\
-(\VJ_h, \Ba_h) + a_2(\VA_h, \Ba_h) + d_2(\Ba_h, r_h)                        &= (\Vg, \Ba_h),\label{full-dis:A}\\
d_2(\VA_h, s_h)    &= 0.\label{full-dis:r}
\end{align}
\end{subequations}
for any $(\varphibf_h, \psi_h, \Ba_h, s_h) \in V_h^2 \times V_h^3 \times V_h^1 \times V_h^0$.

Because $\Div V_h^2 \subseteq V_h^3$ (see \cite{arnold2020}), letting $\psi_h = \Div\VJ_h$ in \eqref{full-dis:Phi}, we have that
\[\NLtwo{\Div\VJ_h} = 0\]
holds, which means divergence-free condition for $\VJ_h$ is satisfied \cite{L2019SIAM}.

Moreover since $\VB_h = \curl\VA_h \in \BH(\Div,\Omega)$, we naturally have $\Div\VB_h = 0$.
Thus the first goal of the present paper is achieved
\begin{equation}\label{eq:div-free}
\Div\VJ_h = 0, \quad \Div\VB_h = 0,\quad \mathrm{in}~~\Omega.
\end{equation}
We point out that $r_h$ is also zero in the domain which can be proven following the lines of continuous variational formulation, thanks to
the condition $\Div\VJ_h = 0$ and $\Div\Vg = 0$.
\vspace{5mm}

\textbf{Remark 1.} The authors in \cite{BAJ2020} state that total current helicity of MHD system vanishes in the steady limit,
which reads
\[\int_\Omega \VJ\cdot\VB = 0\]
Due to the fact $\VA_h \in V_h^1$ and $\curl V_h^1 \subseteq V_h^2$ (see \cite{arnold2020}),
we remark that letting $\varphibf_h = \VB_h = \curl\VA_h$  in \eqref{full-dis:J},  one will have
\begin{equation}
\eta(\VJ_h,\VB_h) - (\phi_h, \Div\VB_h) - (\Vw\times\curl\VA_h,\curl\VA_h) = 0
\end{equation}
which indicates that $\eta(\VJ_h,\VB_h) = 0$ thanks to the equality $\Div\VB_h = 0$.
Therefore the discrete solutions of our present solver also preserve the total magnetic helicity.
\vspace{5mm}

\textbf{Remark 2.} When the velocity $\Vw=0$ in the domain,
MHD kinematics equations turn into the classical eddy current equations.
Thus we also develop a new finite element solver for it such that the discrete current density and discrete magnetic induction are
simultaneously divergence-free.
\section{A block preconditioner}
In this section we will propose a robust block preconditioner based on constraint preconditioning framework for the algebraic systems.
After finite element discretization, we will get the linear algebraic system
\begin{equation}\label{Axb}
\mathcal{A}\Vx = \Vb
\end{equation}
where the vector $\Vx$ consists of the degrees of freedom for $(\VJ_h, \phi_h, \VA_h, r_h)$.
The matrix $\mathcal{A}$ could be written in the following block form
\begin{equation*}
\mathcal{A}=\left(
              \begin{array}{cccc}
              \BM   &G^T      &\BK   &0\\
              G     &0        &0     &0\\
              \BX   &0        &\BF   &B^T\\
              0     &0        &B     &0\\
              \end{array}
            \right)
\end{equation*}
where
\begin{align*}
&\BM_{ij}=\eta(\varphibf_j, \varphibf_i),    \quad\forall \varphibf_i, \varphibf_j \in V_h^2\\
&G^T_{ij}=-(\psi_j, \Div\varphibf_i),        \quad\forall \varphibf_i\in V_h^2, \psi_j \in V_h^3 \\
&\BK_{ij}=(\curl\Ba_j\times\Vw,\varphibf_i), \quad\forall \varphibf_i\in V_h^2, \forall\Ba_j \in V_h^1\\
&\BX_{ij}=-(\varphibf_j, \Ba_i),             \quad\forall \Ba_i      \in V_h^1, \forall\varphibf_j \in V_h^2\\
&\BF_{ij}=\nu_m(\curl\Ba_j, \curl\Ba_i),     \quad\forall \Ba_i,     \Ba_j \in V_h^1\\
&B^T_{ij}=(\nabla s_j, \Ba_i),               \quad\forall \Ba_i      \in V_h^1, \forall  s_j \in V_h^0
\end{align*}
For multi-physics problems, block preconditioning is famous \cite{ESW2014}.
Now we attempt to give an efficient block preconditioner motivated by constraint preconditioning theory \cite{KGW2000, C2002}.
\subsection{Constraint preconditioning}
Denote
\begin{equation}
\bbZ = \left(
  \begin{array}{cc}
  \BM     &\BK \\
  \BX     &\BF \\
  \end{array}
\right),\quad \bbN = \left(
  \begin{array}{cc}
  G     &0 \\
  0     &B \\
  \end{array}
\right)
\end{equation}
and we rearrange the order of sub-matrices as follows for clear explanation
\begin{equation}
\widetilde{\mathcal{A}} = \left(
  \begin{array}{cc}
  \bbZ     &\bbN^\bbT \\
  \bbN     & 0 \\
  \end{array}
\right)
\end{equation}
Let $\widetilde{\bbZ}$ be a approximation of $\bbZ$, we can obtain a block matrix
\begin{equation}\label{eq:CT1}
\widetilde{\mathcal{P}} = \left(
  \begin{array}{cc}
  \widetilde{\bbZ}     &\bbN^\bbT \\
  \bbN     & 0 \\
  \end{array}
\right)
\end{equation}
The basic ideas of constraint preconditioning say that if $\widetilde{\bbZ}$ is a good approximation of $\bbZ$,
then $\widetilde{\mathcal{P}}$ could be a good preconditioner for  $\widetilde{\mathcal{A}}$ \cite{KGW2000, C2002}.
Let the order of $\bbZ$ is $N_F\times N_F$ and the order of $\bbN$ is $N_L\times N_F$. Consider the following
generalized eigenvalue problems,
\begin{equation}
\widetilde{\mathcal{A}}\Vx = \lambda\widetilde{\mathcal{P}}\Vx
\end{equation}
in a more rigorous way, constraint preconditioning says that $\widetilde{\mathcal{P}}^{-1}\widetilde{\mathcal{A}}$
has eigenvalues 1 with multiplicity $2N_L$ and the remaining $N_F - N_L$ eigenvalues are those of
$\bbS = (\bbV^T\widetilde{\bbZ}\bbV)^{-1}(\bbV^T\bbZ\bbV)$. Here $\bbV$ is composed of the
orthogonal basis of null space of matrix $\bbN$ \cite{KGW2000, C2002}.

Thus if eigenvalues of $\bbS$ is clustered tightly, then that is case of
$\widetilde{\mathcal{P}}^{-1}\widetilde{\mathcal{A}}$. Considering the eigenvalue problem $\bbS\Vy = \lambda\Vy$,
we have
\begin{equation}
\Vy^T \bbV^T \bbZ \bbV \Vy = \lambda \Vy^T \bbV^T \widetilde{\bbZ} \bbV \Vy
\end{equation}
Denote by $\Vx = \bbV \Vy$, we equivalently have
\begin{equation}\label{eq:CT2}
\Vx^T\bbZ\Vx = \lambda\Vx^T\widetilde{\bbZ}\Vx,\quad \bbN\Vx = \mathbf{0}, \quad \forall \Vx \in \bbR^{N_F}
\end{equation}
Due to \eqref{eq:CT2}, our main idea is to first find a good approximation of matrix $\bbZ$. Note that the LU decomposition
of $\bbZ$ reads
\[\bbZ = \left(
  \begin{array}{cc}
  \BM     &\BK \\
  \BX     &\BF \\
  \end{array}
\right)=\left(
  \begin{array}{cc}
  \BI     &0 \\
  \BX\BM^{-1}     & \BI \\
  \end{array}
\right)\left(
  \begin{array}{cc}
  \BM     &\BK \\
  0     & \BF-\BX\BM^{-1}\BK \\
  \end{array}
\right)\]
It can be expected that the upper block matrix should be a good preconditioner for $\bbZ$, namely one can choose
\begin{equation}
\widetilde{\bbZ} = \left(
  \begin{array}{cc}
  \BM     &\BK \\
  0     & \BF-\BX\BM^{-1}\BK \\
  \end{array}
\right)
\end{equation}
But the inverse of $\BM^{-1}$ is difficult to use in practice.
We observe that the underlying operator of $\BF-\BX\BM^{-1}\BK$ is
\[\nu_m\curl\curl(\cdot) + \sigma\curl(\cdot)\times\Vw\]
So we choose $\BF_w$ of the following variational formula to approximate $\BF-\BX\BM^{-1}\BK$
\[\BF_w[i,j]:=\nu_m(\curl\Ba_j,\curl\Ba_i) + \sigma(\curl\Ba_j, \Vw\times\Ba_i)\]
Namely we choose $\widetilde{\bbZ}$ to be
\begin{equation}\label{eq:CT3}
\widetilde{\bbZ} = \left(
  \begin{array}{cc}
  \BM     &\BK \\
  0     & \BF_w \\
  \end{array}
\right)
\end{equation}
Due to \eqref{eq:CT3} an initial preconditioner $\mathcal{A}_1$ can be obtained for original matrix $\mathcal{A}$
\begin{equation}\label{mat:A1}
\mathcal{A}_1=\left(
              \begin{array}{cccc}
              \BM &G^T      &\BK   &0\\
              G   &0        &0     &0\\
              0   &0        &\BF_w   &B^T\\
              0   &0        &B     &0\\
              \end{array}
            \right)
\end{equation}
However for \eqref{mat:A1}, we still need to consider the efficient preconditioning for the following two saddle systems
\begin{equation}\label{mat:JA}
\mathcal{A}_J=\left(
  \begin{array}{cc}
  \BM   &G^T \\
  G     &0   \\
  \end{array}
\right),\quad \mathcal{A}_a = \left(
                \begin{array}{cc}
                \BF_w  &B^T \\
                B    &0   \\
                \end{array}
              \right).
\end{equation}
\subsection{Preconditioning for $\mathcal{A}_J$}
In this subsection we will devise efficient preconditioning for $\mathcal{A}_J$, which is the matrix of mixed method for elliptic equations.
From the references \cite{vassi1996, arnold97}, the authors proposed the following block diagonal preconditioner
\begin{equation}\label{pre:JA}
\left(
  \begin{array}{cc}
  \widehat{\BM}  &0 \\
  0           &\sigma Q           \\
  \end{array}
\right)
\end{equation}
where
\[\widehat{\BM}_{ij}:=\sigma^{-1}(\varphibf_j, \varphibf_i) + \sigma^{-1}(\Div\varphibf_j,\Div\varphibf_i),\quad Q_{ij}:=(\psi_j,\psi_i)\]
However, here we will give some further modifications to improve \eqref{pre:JA}.
These improvements are based on augmentation, approximate
block decompositions and the commutativity of the underlying continuous operators.
Note that
\begin{equation}\label{AL:J}
\left(
    \begin{array}{cc}
      I & \eta G^T Q^{-1} \\
      0 & I \\
    \end{array}
  \right)\left(
           \begin{array}{cc}
           \BM & G^T \\
            G & 0 \\
           \end{array}
         \right) = \left(
                     \begin{array}{cc}
                       \BM + \eta G^T Q^{-1} G  & G^T \\
                       G & 0 \\
                     \end{array}
                   \right)
\end{equation}
and
\[\left(
  \begin{array}{cc}
  \BM + \eta G^T Q^{-1} G  & G^T \\
  G & 0 \\
  \end{array}
\right)=\left(
          \begin{array}{cc}
            I & 0 \\
            G\widetilde{\BM}^{-1} & I \\
          \end{array}
        \right)\left(
                 \begin{array}{cc}
                   \widetilde{\BM} & G^T \\
                   0               & - G\widetilde{\BM}^{-1} G^T \\
                 \end{array}
               \right) \triangleq \mathcal{L}_J\mathcal{U}_J
\]
with $\widetilde{\BM} = \BM + \eta G^T Q^{-1} G$.
The underlying continuous operator of $\widetilde{\BM} $ is
\[\eta \mathrm{Id} - \eta \nabla\Div\]
so $\widehat{\BM}$ should be an ideal approximation for $\widetilde{\BM}$.
The operator of $G\widetilde{\BM}^{-1} G^T$ is as follows
\[-\Div(\eta \mathrm{Id}- \eta\nabla\Div)^{-1}\nabla\]
Because Laplace operator can commutate with the gradient operator
\[(\eta \mathrm{Id} - \eta\nabla\Div)\nabla = \eta\nabla(\mathrm{Id} - \Delta)_\phi
\Rightarrow \sigma \nabla(\mathrm{Id} - \Delta)_\phi^{-1} = (\eta \mathrm{Id} - \eta\nabla\Div)^{-1}\nabla\]
we obtain
\[-\Div(\eta \mathrm{Id} - \eta\nabla\Div)^{-1}\nabla = -\sigma\Div\nabla (\mathrm{Id} - \Delta)_\phi^{-1}
= \sigma(-\Delta_\phi)(\mathrm{Id} - \Delta)_\phi^{-1}\]
The operators $-\Delta_\phi$ and $(\mathrm{Id} - \Delta)_\phi$ are spectrally equivalent,
so we can use identity operator $\mathrm{Id}$ to
approximate $(-\Delta_\phi)(\mathrm{Id} - \Delta)_\phi^{-1}$. Thus
\[-\Div(\eta \mathrm{Id} - \eta\nabla\Div)^{-1}\nabla \approx \sigma\mathrm{Id}\]
From above we can approximate $\mathcal{U}_J$ as
\begin{equation}\label{preJ1}
\left(
                 \begin{array}{cc}
                   \widehat{\BM} & G^T \\
                   0 & -\sigma Q \\
                 \end{array}
               \right)
\end{equation}
Due to \eqref{AL:J} and \eqref{preJ1},
\begin{equation}\label{preJ2}
\mathcal{P}_J = \left(
    \begin{array}{cc}
    I & -\eta G^T Q^{-1} \\
    0 & I \\
    \end{array}
  \right)\left(
  \begin{array}{cc}
  \widehat{\BM} & G^T \\
  0 & -\sigma Q \\
  \end{array}
  \right) = \left(
            \begin{array}{cc}
            \widehat{\BM} & 2G^T \\
            0 & -\widehat{Q} \\
            \end{array}
            \right)
\end{equation}
should be a good preconditioner for $\mathcal{A}_J$ with $\widehat{Q} = \sigma Q$.
\subsection{Preconditioning for $\mathcal{A}_a$}
Next we will consider the block preconditioning for $(\VA_h, r_h)$ part.
Let $L$ be the matrix of operator $-\Delta_r$ on $V_h^0$, namely the finite element space for $r_h$.
Note that
\begin{equation}\label{AL:A}
\left(
    \begin{array}{cc}
    I & B^T L^{-1} \\
    0 & I \\
    \end{array}
  \right)\left(
           \begin{array}{cc}
           \BF_w & B^T \\
           B & 0 \\
           \end{array}
         \right) = \left(
                     \begin{array}{cc}
                     \BF_w + B^T L^{-1} B  & B^T \\
                     B & 0 \\
                     \end{array}
                   \right)
\end{equation}
and
\[\left(
  \begin{array}{cc}
  \BF_w + B^T L^{-1} B  & B^T \\
  B & 0 \\
  \end{array}
\right)=\left(
          \begin{array}{cc}
          I & 0 \\
          B\widetilde{\BF}_w^{-1} & I \\
          \end{array}
        \right)\left(
                 \begin{array}{cc}
                 \widetilde{\BF}_w & B^T \\
                 0 & - B\widetilde{\BF}_w^{-1} B^T \\
                 \end{array}
               \right) \triangleq \mathcal{L}_a\mathcal{U}_a
\]
with $\widetilde{\BF}_w = \BF_w + B^T L^{-1} B$. As above we want to give a good matrix approximation for $\mathcal{U}_a$.
In the references \cite{gre07,phil2014},
the authors use mass matrix to approximate $B^T L^{-1} B$. Inspired by their work, we use
the following matrix $\widehat{\BF}_w$
\[\widehat{\BF}_w[i,j] : = \nu_m(\curl\Ba_j,\curl\Ba_i) + \sigma(\curl\Ba_j, \Vw\times\Ba_i) + (\Ba_j,\Ba_i)\]
to approximate $\widetilde{\BF}_w$.

Next consider the continuous Schur complement operator of $B\widetilde{\BF}_w^{-1} B^T$
\[-\Div[\nu_m\curl\curl + \sigma\curl(\cdot)\times\Vw+ \nabla(-\Delta)^{-1}(-\Div)]^{-1}\nabla\]
Note that
\[[\nu_m\curl\curl + \sigma\curl(\cdot)\times\Vw + \nabla(-\Delta)^{-1}(-\Div)]\nabla = \nabla\]
Due to $\curl\nabla = 0$ thus we have
\[[\nu_m\curl\curl + \sigma\curl(\cdot)\times\Vw + \nabla(-\Delta)^{-1}(-\Div)]^{-1}\nabla = \nabla\]
and further
\[-\Div[\nu_m\curl\curl + \sigma\curl(\cdot)\times\Vw + \nabla(-\Delta)^{-1}(-\Div)]^{-1}\nabla = -\Div\nabla = -\Delta_r\]
According to the discussions above, the matrix $L$ should be a good approximation for $B\widetilde{\BF}_w^{-1} B^T$.
Finally a reasonable approximation of $\mathcal{U}_a$ is as follows
\begin{equation}\label{preA1}
\left(
\begin{array}{cc}
\widehat{\BF}_w &B^T \\
0             &-L  \\
\end{array}
\right)
\end{equation}
And due to \eqref{AL:A}
\begin{equation}\label{preA2}
\mathcal{P}_a = \left(
    \begin{array}{cc}
      I & -B^T L^{-1} \\
      0 & I           \\
    \end{array}
  \right)\left(
  \begin{array}{cc}
  \widehat{\BF}_w & B^T \\
  0             & -L \\
  \end{array}
  \right) = \left(
  \begin{array}{cc}
  \widehat{\BF}_w & 2B^T \\
  0             & -L   \\
  \end{array}
  \right)
\end{equation}
should be a good preconditioner for $\mathcal{A}_a$.

\subsection{Preconditioned FGMRES solver}
From \eqref{mat:A1}, \eqref{preJ2}, \eqref{preA2} and the derivation above,
now we present our final block preconditioner for $\mathcal{A}$
\begin{equation}\label{pre1:A}
\mathcal{P} =\left(
              \begin{array}{cccc}
              \widehat{\BM}   &2G^T         &\BK           &0    \\
              0               &-\widehat{Q} &0             &0    \\
              0               &0            &\widehat{\BF}_w &2B^T \\
              0               &0            &0             &-L   \\
              \end{array}
            \right)
\end{equation}
In our preconditioning implementation, iterative solvers are used. So small changes exist in the preconditioning at every step.
And we choose FGMRES \cite{saad1993} to solve \eqref{Axb}.
Given a general vector $\Vx$ which has the same size as one column vector of $\mathcal{A}$,
we let $(\Vx_J, \Vx_\phi, \Vx_A, \Vx_r)$ be the vectors which consist of entries of $\Vx$ corresponding to
$(\VJ_h, \phi_h, \VA_h, r_h)$ respectively. Denote by $\Vr = (\Vr_J, \Vr_\phi, \Vr_A, \Vr_r)$ the residual vector from FGMRES solver,
then in every FGMRES iteration one needs to solve the following preconditioning equations
\[\mathcal{P}\Ve = \Vr\]
where $\Ve = (\Ve_J, \Ve_\phi, \Ve_A, \Ve_r)$. We solve this preconditioning equations by the following approximate iterative methods:
\begin{enumerate}
  \item Solve $L \Ve_r = -\Vr_r$ by preconditioned CG method with relative tolerance $\varepsilon_0$.
  The preconditioner is the algebraic multigrid method (AMG) \cite{hen02}.
  \item Solve $\widehat{\BF}_w\Ve_A = \Vr_A - 2B^T\Ve_r$ by preconditioned GMRES with relative tolerance $\varepsilon_0$.
  The preconditinoer is one-level additive Schwarz method.
  \item Solve $\widehat{Q} \Ve_ \phi = -\Vr_{\phi}$ by 5 CG iterations with the diagonal preconditioning.
  \item Solve $\widehat{\BM} \Ve_J = \Vr_J - 2G^T \Ve_{\phi} - \BK\Ve_A$ .
  We use preconditioned CG method with the auxiliary space preconditioner \cite{hip07} using relative tolerance $\varepsilon_0$.
\end{enumerate}
\section{Numerical experiments}
In this section, we will present two numerical examples to verify the divergence-free feature of the discrete solutions,
the convergence rate of finite element solutions and the performance of the block preconditioner.
Due to $\VB_h = \curl\VA_h \in \BH(\Div,\Omega)$ we naturally have $\NLtwo{\Div\VB_h} = 0$.
So we will only report the divergence of the discrete current density $\VJ_h$.

The code is developed based on the finite element package-\textbf{Parallel Hierarchical Grid (PHG)} \cite{zh09-1,zh09-2}.
The computational domain $\Omega$ is a unit cube $(0,1)^3$.
We use PETSc's FGMRES solver \cite{petsc} and set the relative tolerances by $\varepsilon = 10^{-10}$ for $\mathcal{A}\Vx=\Vb$
and $\varepsilon_0 = 10^{-3}$ for sub-solvers in the preconditioner $\mathcal{P}$.

\begin{example}[Precision test]\label{ex2}
In this test example, we choose the velocity field $\Vw=(x,y,z)$ and set $\sigma = {\sf Rm} = 1$.
The following analytic solutions are used to test the convergence rate
\[\VJ = (\sin y, 0, x^2),~\phi = z,~\VA = (0, \cos x , 0),~r=0.\]
\end{example}

\begin{table}[!htb]
  \centering
  \caption{Convergence rate of the finite element solutions (Example \ref{ex2}).}
  \label{EnergyError}
  \begin{tabular*}{14cm}{@{\extracolsep{\fill}}|c|c|c|c|c|c|c|}
  \hline
  $h$    &$\|\VJ-\VJ_h\|_{\BH(\Div)}$  &order &$\|\phi - \phi_h\|_{L^2}$  &order &$\|\VA-\VA_h\|_{\BH(\curl)}$ &order\\ \hline
  0.86603 &5.9811e-02   &---           &1.0208e-01   &---     &9.8060e-02  &---     \\\hline
  0.43301 &2.6438e-02   &1.1778        &5.1034e-02   &1.0002  &4.8104e-02  &1.0275  \\\hline
  0.21651 &1.2527e-02   &1.0776        &2.5516e-02   &1.0001  &2.3780e-02  &1.0164  \\\hline
  0.10825 &6.1235e-03   &1.0326        &1.2758e-02   &1.0000  &1.1821e-02  &1.0084  \\\hline
  0.05413 &3.0326e-03   &1.0138        &6.3789e-03   &1.0000  &5.8932e-03  &1.0042  \\\hline
  \end{tabular*}
\end{table}

\begin{table}[!htb]
  \centering
  \caption{$L^2$ norm of divergence of $\VJ_h$ (Example \ref{ex2}).}
  \label{L2Error}
  \begin{tabular*}{14cm}{@{\extracolsep{\fill}}|c|c|c|c|c|c|}
  \hline
  $h$                     &0.86603      &0.43301     &0.21651    &0.10825    &0.05413  \\\hline
  $\|\Div\VJ_h\|_{L^2}$   &4.0087e-12   &2.5019e-12  &1.6498e-12 &4.7832e-11 &8.2177e-12        \\\hline
  \end{tabular*}
\end{table}

From Table \ref{EnergyError} and Table \ref{L2Error}, we can see that optimal convergence rates are obtained.
And the $\|\Div\VJ_h\|_{L^2}$  is very small compared with the finite element error.
This example shows that the divergence-free conditions for $\VJ_h$ and $\VB_h$ are both satisfied.
The sources of deviation from precisely zero mainly come from the
tolerance $\varepsilon = 10^{-10}$ of solving $\mathcal{A}\Vx = \Vb$.

\begin{example}[Performance of the block preconditioner]\label{ex3}
In this example, we set $\sigma = 1$ and prescribe the velocity field by
\begin{equation}
\Vw = \left(
        \begin{array}{c}
          -16x(1-x)y(1-y)\sin \theta \\
          16x(1-x)y(1-y)\cos \theta \\
          0 \\
        \end{array}
      \right)
\end{equation}
where $\theta$ is the angle between vector $(x,y,0)$ and the positive direction of $x$-axis.
The applied magnetic field $\VB_s = (1,0,0)$ and the boundary condition for $\VA$ is $\VA_b = (0, 0, y)$
such that $\curl\VA_b = \VB_s$. Zero boundary condition for $\phi$ is used. The source term $\Vg$ is zero.
\end{example}

In this example, we want to test the performance of the proposed block preconditioner \eqref{pre1:A}.
The information of grids and degree of freedoms are listed in Table \ref{pre:grid}.
We use three different magnetic Reynolds number $\textsf{Rm} = 50, 100, 200$
to show the performance of the preconditioenr $\mathcal{P}$.

From Table \ref{pre1:tau1}, we observe that the quasi-optimality of $\mathcal{P}$  with respect to the grid refinement.
And the preconditioned FGMRES solver is still robust for relatively high physical parameter $\textsf{Rm} = 200$.
Note that the relative error tolerance for FGMRES is $10^{-10}$.

\begin{table}[!htb]
  \centering
  \caption{The mesh sizes and the numbers of DOFs (Example \ref{ex3}).}
  \label{pre:grid}
  \begin{tabular*}{10cm}{@{\extracolsep{\fill}}|c|c|c|c|}
  \hline
  Mesh        &$h$         &DOFs for $\VJ+\phi$       &DOFs for $\VA+r$  \\\hline
  $\Ct_1$     &0.86603     &360+48                    &196+125                 \\\hline
  $\Ct_2$     &0.43301     &2592+384                  &1208+729                \\\hline
  $\Ct_3$     &0.21651     &19584+3072                &8368+4913               \\\hline
  $\Ct_4$     &0.10825     &152064+24576              &62048+35937               \\\hline
  $\Ct_5$     &0.05413     &1198080+196608            &477376+274625               \\\hline
  \end{tabular*}
\end{table}
\begin{table}[!htb]
  \centering
  \caption{FGMRES iteration number with preconditioner $\mathcal{P}$ (Example \ref{ex3}).}
  \label{pre1:tau1}
  \begin{tabular*}{8cm}{@{\extracolsep{\fill}}|c|c|c|c|}
  \hline
  Mesh        &\textsf{Rm} = 50   &\textsf{Rm} = 100  &\textsf{Rm} = 200\\\hline
  $\Ct_1$     &21      &23    &30   \\\hline
  $\Ct_2$     &19      &22    &30   \\\hline
  $\Ct_3$     &16      &18    &25   \\\hline
  $\Ct_4$     &14      &16    &20   \\\hline
  $\Ct_5$     &14      &15    &19   \\\hline
  \end{tabular*}
\end{table}
\section{Conclusions}
In this paper, we develop a new conservative finite element iterative solver
for the three-dimensional MHD kinematics equations which can ensure
exactly divergence-free approximations of the current density and magnetic induction.
The magnetic helicity is also preserved by the solver.
Moreover, we devise a robust block preconditioner motivated by the constraint preconditioning framework \cite{KGW2000, C2002}.
Future investigation will consider the efficient implicit solver based on the proposed algorithms.
Moreover, combined with our previous
work in \cite{L2019SIAM} the implementation of full MHD solver using $(\vec{u},p,\vec{J},\phi,\vec{A},r)$ as solution variables is
under consideration. Our initial results shows that the ideas of preconditioning of the present work still seems
to be applicable for the new systems.

\section*{Acknowledgments}
\textbf{Lingxiao Li} was supported by National Natural Science Foundation of China under Grant 11901042.
\textbf{Xujing Li} was supported in part by National Natural Science Foundation of China under Grant 11771440 and 11805049.

\end{document}